\documentclass[10pt,twocolumn,twoside]{IEEEtran}
\usepackage{graphicx}
\usepackage{epsfig}

\usepackage{amssymb,amsmath,amsthm, amsbsy}
\usepackage[UKenglish]{babel}

\usepackage{algorithm} 
\usepackage{algpseudocode}

\usepackage{mathtools}
\usepackage{array}

\usepackage{changepage} 

\usepackage{enumitem}

\usepackage{graphicx}
\graphicspath{{figuresPartI/}} 
\usepackage[export]{adjustbox}

\usepackage{siunitx}

\usepackage{cases} 
\usepackage[usenames, dvipsnames]{color}

\newcommand{\scalevariableS}{0.73}

\usepackage{xcolor}

\definecolor{line_gray}{gray}{0.75}
\definecolor{cell_gray}{gray}{0.9}
\definecolor{cell_grad_blue1}{HTML}{42A5F5}
\definecolor{cell_grad_blue2}{HTML}{64B5F6}
\definecolor{cell_grad_blue3}{HTML}{90CAF9}
\definecolor{cell_grad_blue4}{HTML}{BBDEFB}

\usepackage{colortbl}

\usepackage{multirow}

\newtheorem{theorem}{Theorem}
\newtheorem{lemma}{Lemma}
\newtheorem{definition}{Definition}
\newtheorem{remark}{Remark}
\newtheorem{assumption}{Assumption}
\newtheorem{requirement}{Requirement}
\newtheorem{corollary}{Corollary}

\newtheorem{constraints}{Constraints}

\usepackage{booktabs}

\newcommand{\lsym}{L_{sym}}

\newcommand{\mati}{\Phi}
\newcommand{\matii}{\Psi}
\newcommand{\espacei}{\mathcal{W}}
\newcommand{\espaceii}{\mathcal{V}}
\newcommand{\evecsi}{W}
\newcommand{\evecsii}{V}
\newcommand{\evali}{\phi}
\newcommand{\evalii}{\psi}

\newcommand{\dimevec}{r}

\def\EqualDef{\mathrel{\mathop=^{\rm def}}}

\newcommand\numberthis{\addtocounter{equation}{1}\tag{\theequation}}

\newcommand{\eqnref}[1]{(\ref{#1})}

\setlist[enumerate]{label*=\arabic*.} 

\begin{document}

\title{Extending the Davis--Kahan theorem for comparing eigenvectors of two symmetric matrices I: Theory}

\author{J.~F.~Lutzeyer and A.~T.~Walden,~\IEEEmembership{Senior~Member,~IEEE}  
\thanks{
Copyright (c) 2019 IEEE. Personal use of this material is permitted. However, permission to use this material for
any other purposes must be obtained from the IEEE by sending a request to pubs-permissions@ieee.org. 
J.~F.~Lutzeyer and A.~T.~Walden are with the Dept. of Mathematics, Imperial College London, London SW7 2AZ, UK (e-mail: jl7511@imperial.ac.uk and a.walden@imperial.ac.uk)
} }
\IEEEpubid{}
\maketitle

\begin{abstract}
The Davis--Kahan theorem can be used to bound the distance of the spaces spanned by the first $r$ eigenvectors of any two symmetric matrices.  
We extend the Davis--Kahan theorem to apply to the comparison of the union of eigenspaces of any two symmetric matrices by making use of polynomial matrix transforms 
and in so doing, tighten the bound. The transform allows us to move requirements present in the original Davis--Kahan theorem, from the \textit{eigenvalues} of the compared matrices on to the \textit{transformation parameters}, with the latter being under our control. 
We provide a proof of concept example, comparing the spaces spanned by the unnormalised and normalised graph Laplacian eigenvectors for $d$-regular graphs, in which the correct transform is automatically identified. 
\end{abstract}

\section{Introduction}

The Davis--Kahan (DK) theorem is a tool for comparing the spaces spanned by the eigenvectors of two symmetric matrices, given that the corresponding eigenvalues satisfy a certain structure.  In this paper, we observe that use can be made of matrix transformations  to extend the applicability of the Davis--Kahan theorem by removing the straitjacket of such a highly constrained eigenvalue structure. The value of the transformation step can be seen when thinking about the many practical uses for such a comparison. Graphs can be represented in multiple ways via different graph shift operator matrices, for example the adjacency matrix, and normalized and unnormalized Laplacians. Bounding the distance of the subspaces spanned by the eigenvectors of these graph shift operators throws up some interesting issues, e.g., the largest eigenvalues of the adjacency matrix correspond to the smallest eigenvalues of the Laplacians, and an eigenvector comparison corresponds to comparing opposite ends of the eigenvalue spectrum; this renders the standard DK theorem inapplicable (due to the eigenvalue restrictions), while our new extended version incorporating a matrix transformation is highly apposite. Comparing spaces spanned by the eigenvectors of the graph shift operators is of particular interest to the signal processing community, since these eigenvectors give rise to the different, much utilised, graph Fourier transforms (GFT) \cite{Shuman2013, Dong2019}.
The GFT is utilised by \cite{Wang2018} for complex brain network analysis, by \cite{Singh2016} to understand centrality patterns over several types of networks and \cite{Lagunas2018} define a graph similarity distance from the graph Fourier bases. 
As stated in \cite[p.~816]{Ortega2018}, making the appropriate graph shift operator choice for a given application remains an open issue.
The work in the current paper provides an improved analytical understanding of the difference between the graph shift operators.

A second related application is the comparison of spaces spanned by eigenvectors of the graph shift operators to spaces spanned by the eigenvectors of their generating matrices in the case of the ubiquitous stochastic block model; this can be used for consistency and rate of convergence studies for different methods based on the eigenvectors of the graph shift operators in such a model  \cite{Cape2019, Eldridge2018, Rohe2011}. 
A third example is the comparison of spaces spanned by the eigenvectors of the sample covariance matrix and its population covariance matrix in a spiked covariance model \cite{Johnstone2001, Jung2009}. In these second and third examples our new approach can deliver tighter bounds than the standard DK theorem, due to the matrix transform allowing a relaxation and utilisation of the eigenvalue structure imposed by the standard DK theorem.

To describe the novel contributions of  this work, we need to introduce our standard notation.
We denote the two symmetric matrices of interest by $\mati$ and $\matii.$ 
The matrices formed from the eigenvectors of these matrices will be denoted by $\evecsi$ and $\evecsii,$ respectively. 
We work with two commonly considered unitarily invariant distance metrics to measure the distance of subspaces of $\mathbb{R}^n$. The considered metrics are $\| (I - \evecsi \evecsi^T) \evecsii \evecsii^T)\|_2$ and $\inf_{R \in O(\dimevec)} \| \evecsi - \evecsii R\|_F,$
where $O(\dimevec)$ denotes the orthogonal group of $\dimevec \times \dimevec $ orthogonal matrices.  Both of these metrics, which use the two norm or spectral norm and Frobenius norm respectively, (defined later),  are functionally related to the canonical angles between spaces and later we shall relate them via an inequality. 

The standard DK theorem  \cite[p.~211-2]{Bhatia2013} enables bounding of the metric $\| (I - \evecsi \evecsi^T) \evecsii \evecsii^T)\|_2,$ measuring the distance between spaces spanned by eigenvectors of the matrices, given that their corresponding eigenvalues satisfy the specified standard structure.
By applying a polynomial transformation to one of the matrices under comparison we leave $\| (I - \evecsi \evecsi^T) \evecsii \evecsii^T)\|_2$ unchanged, since the eigenvectors are unchanged, but can change the location of the corresponding eigenvalues.
This allows us to relax the necessary DK structure of the eigenvalues 
to a non-zero eigengap structure.
Furthermore, we find that since the DK bound 
depends on the matrices and eigenvalues under comparison, we can choose the polynomial transformation to not only relax the conditions on the eigenvalues, but also to reduce the DK bound value. 
Use of the metric $\inf_{R \in O(\dimevec)} \| \evecsi - \evecsii R\|_F$ is common \cite{Lei2015,Rohe2011,Tang2018} and 
we are able to produce an upper bound for it  using our sharpened version of the DK theorem.

In Section~\ref{sec:background} we briefly discuss previous work on the DK theorem. In Section~\ref{sec:subspaces} we give an interpretation of two distance measures on subspaces spanned by eigenvectors, and relate them via an inequality. In Section
\ref{sec_proof_structure} we show that under a suitable assumption, we are always able to compare the spaces spanned by the first $r$ eigenvectors of two symmetric matrices via
the DK theorem; the bound involves a between-matrix eigengap. Section~\ref{sec_eval_map} presents an extension to the DK theorem which uses a polynomial matrix transformation of one of the symmetric matrices and also  applies to
the comparison of any 
two sets of consecutive and corresponding  eigenvectors from two symmetric matrices. 
Assumptions on the spectra are developed. The special case of affine transforms is discussed in Section~\ref{sec:affinetransforms} and a proof of concept is given in Section~\ref{sec_d_reg_example}. Our summary and conclusions are given in Section~\ref{sec:sum}.

In paper II we tackle the computational issues with a focus on affine transforms where
a fractional programming approach is used. Such problems can be transformed to convex optimisation problems 
with a unique global solution.
We calculate our extended DK bound in situations
where the standard DK eigenvalue structure is not satisfied and the standard DK theorem is inapplicable, (the graph shift operators example mentioned above), but also introduce two situations in which standard DK bounds are available 
but our bound produced via the transformation approach 
attains lower bound values (the generating matrices and covariance matrices examples mentioned above).

\section{Some Background on the DK theorem}\label{sec:background}

The DK theorem first appeared in 1970  \cite{Davis1970} and is a topic of current research, being of interest in, e.g., the analysis of spectral graph embedding methods \cite{Tang2018} and of principal component analysis of covariance matrices \cite{Fan2013, Wang2017}. 
 von Luxburg \cite{vonLuxburg2007} states that the DK theorem forms the basis for the perturbation approach to spectral clustering of networks, where cluster structure is seen as a perturbed version of a structure with perfectly disjoint clusters.
A new variant with a statistical flavor was recently published by \cite{Yu2015}.

Tightening the DK theorem using different norms, probabilistic methods and making mild assumptions on the structure of the underlying matrices has also been the subject of much recent research work. 
 \cite{ORourke2018} use a probabilistic approach to sharpen the DK theorem using certain structural properties of one of the matrices, such as having low rank. 
In \cite{Cape2019} the DK theorem is used as a coarse benchmark for their bound on the difference of spaces spanned by eigenvectors in the 2-to-infinity norm.
\cite{Eldridge2018} also state that the tightness of the DK theorem is in many settings suboptimal and hence present their own tightened bound using the infinity norm. 
In contrast to these approaches, our generalized and tightened DK theorem applies to any pair of symmetric matrices. 

As well as the tightness issues, a drawback of the DK theorem is that the choice of eigenvalue intervals might not be suitable for the two matrices we want to compare. For example in the particular context of one matrix being considered a perturbation of the other,   \cite[p.~407]{vonLuxburg2007}  von Luxburg noted that ``If the perturbation is too large or the eigengap is too small, we might not find a set $S_1$ such that both the first $r$ eigenvalues of the perturbed and unperturbed versions  ... are contained in $S_1.$''
Indeed, if we want to compare eigenvectors corresponding to eigenvalues lying on opposing ends of the spectra, as for example when comparing corresponding adjacency to graph Laplacian eigenvectors, then, finding eigenvalue intervals such that the DK theorem can be applied is rarely possible. 
By applying a polynomial matrix transform to one of the matrices under comparison, we are able to show that for any two symmetric matrices, given nonzero eigengaps on either side of two contiguous sets of eigenvalues, we are able to find a valid set of eigenvalue intervals such that their corresponding eigenvectors can be bounded using a DK bound.  

To summarize, the  key innovation in this paper is that the application of a matrix transformation allows us to broaden the class of cases in which the DK theorem applies by shifting the restrictive assumptions from the \textit{eigenvalues} 
to the \textit{transformation parameters}, with  
the latter being under our control. As a consequence the DK theorem can then be applied to the comparison of any 
two sets of consecutive and corresponding  eigenvectors from two symmetric matrices (excluding degenerate examples).

\section{Subspaces of $\mathbb{R}^n$ --- spaces spanned by $\dimevec\leq n$ eigenvectors }\label{sec:subspaces}

\subsection{Groundwork}
We shall compare spaces spanned by eigenvectors of two symmetric matrices $\mati, \matii \in \mathbb{R}^{n \times n}$ with eigenvalues $\evali_1 \leq \evali_2 \leq \ldots \leq \evali_n$ and $\evalii_1 \leq \evalii_2 \leq \ldots \leq \evalii_n$ and corresponding eigenvectors $\{ w_1, w_2, \ldots, w_n\}$ and $\{ v_1, v_2, \ldots, v_n\}$, respectively.

 In fact we have reason to compare the subspaces of $\mathbb{R}^n$ spanned by
 $\dimevec$ consecutive eigenvectors of the two symmetric matrices. 
 We define matrices $W_j  = \left[ w_{j+1}, \ldots, w_{j+\dimevec}\right]$ and $V_j  = \left[ v_{j+1}, \ldots, v_{j+\dimevec}\right]$ holding $\dimevec$ consecutive eigenvectors of $\mati$ and $ \matii$, respectively. Note that $\evecsi_0$ and $\evecsii_0$ correspond to the {\it first}\/ $\dimevec$ eigenvectors, i.e., the $r$ eigenvectors corresponding to the $r$ smallest eigenvalues.  The value of $j\in \{0,1, \ldots, n-r\}$ determines a particular contiguous block of $\dimevec$ eigenvectors, i.e, $j$ is an offset or shift  parameter; if the value of $j$ is irrelevant in a particular discussion then the parameter will be suppressed for brevity.

We now formally introduce the different spaces associated with matrices composed of orthogonal columns, which will be utilised when comparing the eigenvectors of symmetric matrices. \cite{Edelman1998} give a very nice introduction to these spaces. Below we summarise their definitions. 

\begin{definition} \label{defn_spaces}
Let $\dimevec \leq n$. We will be working with the following three matrix groups:
\begin{enumerate}
\item The \textit{orthogonal group}, denoted $O(r)$, of $r\times r $ orthogonal matrices.
\item The \textit{Stiefel manifold}, denoted $\mathbb{V}_{n,\dimevec}$, consisting of $n \times \dimevec$ matrices with orthonormal columns.
\item The \textit{Grassmann manifold}, denoted $\mathbb{G}_{n,\dimevec}$, consisting of $\dimevec$-dimensional subspaces of $\mathbb{R}^n$. Elements of the Grassmann manifold are equivalence classes of elements of the Stiefel manifold, where two elements of the Stiefel manifold are equivalent if their columns span the same subspace.
\hfill $\lhd$
\end{enumerate}
\end{definition}

The following Lemma helps us understand the equivalence classes formed on the Stiefel manifold to produce elements on the Grassmann manifold and is furthermore used to motivate our distance metric on eigenvector sets.

\begin{lemma} \label{lemma_existence_of_Q}
If the columns of two elements $\evecsi, \evecsii \in \mathbb{V}_{n,\dimevec}$ span the same subspace of $\mathbb{R}^n$, then, there exists an orthogonal matrix $Q \in O(\dimevec)$ such that $\evecsi = \evecsii Q$.
\end{lemma}
\begin{IEEEproof}
See Appendix~\ref{app_existence_of_Q}. 
\end{IEEEproof}
Note that two elements $\evecsi, \evecsii \in \mathbb{V}_{n,\dimevec}$, the columns of which span the same subspace of $\mathbb{R}^n$, correspond to a single element in the Grassmann manifold.

\begin{definition} \label{rmk_on_orthogonal_projectors}
The orthogonal projector onto the subspace spanned by the eigenvectors $w_{j+1}, \ldots, w_{j+\dimevec}$ can be expressed in terms of $\evecsi_j=[w_{j+1}, \ldots, w_{j+\dimevec}]$ as $\evecsi_j \evecsi_j^T$ and the orthogonal projector onto the complementary subspace is equal to $(I - \evecsi_j \evecsi_j^T)$. \cite[p.~430]{Meyer2000} \hfill $\lhd$
\end{definition}

It is pointed out in \cite[p.~319]{Edelman1998} that each element in the Grassmann manifold, $\mathbb{G}_{n,\dimevec}$ has a corresponding unique orthogonal projector onto its space of the form $VV^T,$ where $V\in \mathbb{V}_{n,\dimevec}$.

\subsection{Measuring Distance Between Subspaces of $\mathbb{R}^n$}

In Definitions \ref{defn_eigenvector_dist_measure} and \ref{defn_eigenvector_dist_measure2} which follow
we define two metrics on the spaces spanned by two eigenvector sets.
Then we introduce the notion of canonical angles between spaces in Definition \ref{defn_canonical_angles} and relate both metrics to the canonical angles in Theorem \ref{thm_stewart_sun_exact_calculation} and Remark \ref{rmk_rho2_to_canonical_angles}. We end by describing the relation of the two metrics in Lemma \ref{lemma_relating_Grassmann_distance_to_DK}.

\begin{definition} \label{defn_eigenvector_dist_measure} \cite[p.~95]{Stewart1990}
Let $\evecsi, \evecsii \in \mathbb{V}_{n, \dimevec}$ be matrices with orthonormal columns which span $\dimevec$ dimensional subspaces  of $\mathbb{R}^n$, $\espacei, \espaceii \in \mathbb{G}_{n,\dimevec}$, respectively. 
We define the unitarily invariant metric on the distance of $\espacei$ and $\espaceii$, denoted $\rho_1(\espacei, \espaceii)$, as,
\begin{equation} \label{eqn_defn_subspace_distance}
\rho_1(\espacei, \espaceii) = \inf_{R \in O(\dimevec)} \|\evecsi -\evecsii R \|_F,
\end{equation}
where $\| \cdot \|_F$ denotes the unitarily invariant Frobenius norm, the square root of the sum of the squared elements.
\hfill $\lhd$
\end{definition} 

In \cite[p.~95]{Stewart1990} $ \rho_1(\espacei, \espaceii)$ is motivated as a metric for the distance between subspaces of $\mathbb{R}^n$ by pointing out that $ \rho_1(\espacei, \espaceii) = 0$ when $\espacei = \espaceii$. This follows directly from Lemma \ref{lemma_existence_of_Q} by taking $R=Q.$

The second metric on subspaces  
is now defined. 

\begin{definition} \label{defn_eigenvector_dist_measure2}
\cite[p.~94]{Stewart1990}
Let $\evecsi, \evecsii \in \mathbb{V}_{n, \dimevec}$ be matrices with orthonormal columns which span $\dimevec$ dimensional subspaces  of $\mathbb{R}^n$, $\espacei, \espaceii \in \mathbb{G}_{n,\dimevec}$, respectively. 
We define the unitarily invariant metric on the distance of $\espacei$ and $\espaceii$, denoted $\rho_2(\espacei, \espaceii)$, as,
$$
\rho_2(\espacei, \espaceii) = \left\| \evecsi \evecsi^T \left(I - \evecsii \evecsii^T \right) \right\|_2,
$$
where $\|\cdot \|_2$ denotes the unitarily invariant two norm or spectral norm, the largest singular value of the matrix, which for any real symmetric matrix corresponds to the maximum of the absolute value of the largest eigenvalue and the absolute value of the smallest eigenvalue.
\hfill $\lhd$
\end{definition}

Recall from Definition \ref{defn_spaces} that each different basis of a subspace corresponds to a Stiefel manifold element, while on the Grassmann manifold all bases of a subspace are represented by a single element.
Note that the metric $\rho_1(\espacei, \espaceii)$ is directly dependent on Stiefel manifold elements $\evecsi, \evecsii,$ while the metric  $\rho_2(\espacei, \espaceii)$ is working with the projectors $\evecsi \evecsi^T, \left(I - \evecsii \evecsii^T \right).$ Since projectors are unique to their corresponding subspaces, they correspond to Grassmann manifold elements. 
In agreement with \cite{Deri2017},
we want to compare Grassmann manifold elements, i.e., entire subspaces rather than just individual bases of these spaces. The metric $\rho_1(\espacei, \espaceii)$ achieves this subspace comparison via the infimum over the orthogonal matrices ensuring that we are considering all bases of the spaces we want to compare. Therefore, both metrics are comparing Grassmann manifold elements. 
In \cite{Cape2019} the minimisation over all orthogonal matrices in \eqnref{eqn_defn_subspace_distance} is described as enabling \textit{basis alignment} of the two spaces under comparison.

Note that in an orthogonally invariant norm such as the Frobenius norm the direct comparison of individual bases is not possible. 
The orthogonal invariance of the norm implies that the comparison of two elements on the Stiefel manifold corresponds to the comparison of a class of elements of the Stiefel manifold related via orthogonal transformations, i.e., $\left\| \evecsi - \evecsii \right\|_F = \left\| \evecsi Q' - \evecsii Q' \right\|_F$ for all $Q' \in O(\dimevec)$.  Hence, without the minimisation over all orthogonal matrices in  \eqnref{eqn_defn_subspace_distance}  we are simply considering transformations of both $\evecsi$ and $\evecsii$ by the same orthogonal transforms rather than by decoupled orthogonal transformations as achieved by the addition of the minimisation step.

When discussing the distance of subspaces of $\mathbb{R}^n$ the notion of canonical angles, as generalisation of angles between lines, is integral. We will therefore define canonical angles in Definition \ref{defn_canonical_angles} and then discuss how both $\rho_1(\espacei, \espaceii)$ and $\rho_2(\espacei, \espaceii)$ are functionally related to the canonical angles. The following definition is adapted from \cite{Vu2013}.

\begin{definition} \label{defn_canonical_angles}
Let $\espacei, \espaceii \in \mathbb{G}_{n,\dimevec}$ be $\dimevec$ dimensional subspaces of $\mathbb{R}^n$ with orthogonal projectors $\evecsi \evecsi^T$ and $\evecsii \evecsii^T.$ Denote the singular values of $\evecsi \evecsi^T (I - \evecsii \evecsii^T)$ by $\beta_1 \geq \ldots \geq \beta_n$. The \textit{canonical angles} between $\espacei$ and $\espaceii$ are the numbers,
$$
\theta_k(\espacei, \espaceii) = \arcsin(\beta_k).
$$
In the literature the diagonal matrix $\Theta(\espacei, \espaceii) = \mathrm{diag}(\beta_1, \ldots, \beta_n)$ is often considered. 
\hfill $\lhd$
\end{definition}

The following theorem states the functional relationship of $\rho_1(\espacei, \espaceii)$ to the canonical angles between $\espacei$ and $\espaceii.$ 

\begin{theorem} \label{thm_stewart_sun_exact_calculation}
\cite[p.~95]{Stewart1990} 
If $\alpha_i$ is the cosine of the $i^{\rm th}$ canonical angle between $\espacei$ and $\espaceii$, then, 
\begin{equation} \label{eqn_stewart_sun_exact_calculation}
\rho_1(\espacei, \espaceii) = 
\left[2 \sum_{i=1}^n \left( 1- \alpha_i\right)\right]^{1/2}. 
\end{equation}
\hfill $\lhd$
\end{theorem}

Theorem \ref{thm_stewart_sun_exact_calculation} will be of great help when calculating distances between eigenspaces as it provides an exact formula, which allows us to avoid the minimisation over all unitary matrices. The $\alpha_i$'s in Theorem \ref{thm_stewart_sun_exact_calculation} can be calculated via Definition \ref{defn_canonical_angles}. Alternatively, \cite[p.~45]{Stewart1990} state that the cosines of the canonical angles, i.e., the $\alpha_i$'s, between $\espacei, \espaceii \in \mathbb{G}_{n,\dimevec}$ are equal to the singular values of $\evecsii^T \evecsi$. Since $\evecsii^T \evecsi$ is only a $\dimevec \times \dimevec$ matrix, it is preferable to obtain the $\alpha_i$'s from $\evecsii^T \evecsi$, instead of the $n\times n$ projector $\evecsi \evecsi^T (I - \evecsii \evecsii^T)$ from Definition \ref{defn_canonical_angles}.

\begin{remark} \label{rmk_rho2_to_canonical_angles}
Canonical angles are defined via the singular values of the projector $\evecsi \evecsi^T \left(I - \evecsii \evecsii^T \right)$ in Definition \ref{defn_canonical_angles}. Hence, it trivially follows from the definition of the two norm that $\rho_2(\espacei, \espaceii)$ is equal to the sin of the largest canonical angle between $\espacei$ and $\espaceii.$
\hfill $\lhd$
\end{remark}

In Lemma \ref{lemma_relating_Grassmann_distance_to_DK} we relate the two distance metrics $\rho_1(\espacei, \espaceii)$ and $\rho_2(\espacei, \espaceii).$ 
Our proof of Lemma \ref{lemma_relating_Grassmann_distance_to_DK} below agrees with the first steps of the proof of Lemma 5.1 in \cite[p.~232]{Lei2015} except for the addition of the $\min(\dimevec, n-\dimevec)$ term instead of $\dimevec$, since we do not want to exclude cases where $\dimevec > n-\dimevec$. 

\begin{lemma} \label{lemma_relating_Grassmann_distance_to_DK}
Take
$c_{n,r} \,\displaystyle{\EqualDef} \, \surd[2 \min(\dimevec, n-\dimevec)]$ 
and let $\evecsi, \evecsii \in \mathbb{V}_{n,\dimevec}.$ 
There exists a $Q \in O(\dimevec)$ such that, 
\begin{equation} \label{eqn_cost_structurea}
\left\|\evecsi - \evecsii Q\right\|_F\leq \!c_{n,r} \left\|\evecsi \evecsi^T (I - \evecsii \evecsii^T) \right\|_2.
\end{equation}
\end{lemma}
\begin{IEEEproof}
See Appendix~\ref{app:relatemetrics}.
\end{IEEEproof}

\section{Application of the DK theorem to comparison of spaces spanned by the first $\dimevec$ eigenvectors } \label{sec_proof_structure}
We firstly formally introduce the DK theorem in Section~\ref{sec_DK_statement} and then in Section \ref{sec_bounding_first_r_evecs}, we  discuss its application to the comparison of the spaces spanned by the first $\dimevec$ eigenvectors of two matrices and show that given non-zero $\dimevec^{\mathrm{th}}$ eigengaps in both spectra, this comparison can always be made. 

\subsection{The Davis--Kahan Theorem} \label{sec_DK_statement}

The DK theorem \cite{Davis1970} in the form given in \cite{Bhatia2013} is reexpressed in Theorem \ref{thm_davis_kahan_Bhatia}, where we also give the orthogonal projections in terms of the matrix eigenvectors. 
\begin{theorem} \label{thm_davis_kahan_Bhatia} \textit{Davis--Kahan Theorem} \cite[p.~211--212]{Bhatia2013}.
Let $\mati, \matii \in \mathbb{R}^{n \times n}$ be symmetric matrices. Let $S_1$ be an interval  $[a,b]$ and $S_2$ be the complement in $\mathbb{R}$ of the interval $(a- \delta, b+\delta)$, i.e., the intervals $S_1$ and $S_2$ lie a distance $\delta>0$ apart. 
Let the columns of matrix $\evecsi$ be orthonormal eigenvectors corresponding to the eigenvalues of $\mati$ contained in $S_1$ and $\evecsii$ have its columns made up of orthonormal eigenvectors corresponding to the eigenvalues of $\matii$ \textit{not} contained in $S_2$.
Then, for every unitarily invariant norm, denoted $\|\cdot\|$,
\begin{equation} \label{eqn_Bhatia_dk}
\left\| \evecsi ~ \evecsi^T (I - \evecsii~ \evecsii^T) \right\|
\leq \frac{1}{\delta} \left\| \mati - \matii \right\|.
\end{equation}
\hfill $\lhd$
\end{theorem}
Note that the interval definitions are in terms of the parameter triplet $(a, b, \delta).$
From Theorem \ref{thm_davis_kahan_Bhatia} we learn that for the comparison of the spaces spanned by eigenvectors $[ w_{j+1}, \ldots, w_{j+\dimevec}]$ of $\mati$ to $[ v_{j+1},  \ldots, v_{j+\dimevec}]$ of $\matii$ we require the following conditions on the corresponding eigenvalues. 
\begin{requirement} \label{rmk_dk_eigenvalue_requirement}
For the DK theorem to apply to the comparison of the eigenvector matrices $W_j, V_j \in \mathbb{V}_{n,\dimevec}$, we need to choose DK intervals, $S_1 = [a,b], S_2 = \mathbb{R}\backslash (a - \delta, b+\delta)$ for some interval separation $\delta>0$, such that, either 
\begin{align*}
\evali_{j+1}, \ldots, \evali_{j+\dimevec} &\in S_1; \\
\quad \evali_{1}, \ldots, \evali_j, \evali_{j+\dimevec+1}, \ldots, \evali_{n} &\notin S_1 ;  \\
\quad \evalii_{j+1}, \ldots, \evalii_{j+\dimevec} &\notin S_2;\\
 \quad  \evalii_{1}, \ldots, \evalii_j, \evalii_{j+\dimevec+1}, \ldots, \evalii_{n} &\in S_2
\end{align*}
or, by swapping $S_1$ and $S_2$, 
\begin{align*}
\evalii_{j+1}, \ldots, \evalii_{j+\dimevec} &\in S_1; \\
\quad \evalii_{1}, \ldots, \evalii_j, \evalii_{j+\dimevec+1}, \ldots, \evalii_{n} &\notin S_1 ; \\
\quad \evali_{j+1}, \ldots, \evali_{j+\dimevec} &\notin S_2; \\
\quad  \evali_{1}, \ldots, \evali_j, \evali_{j+\dimevec+1}, \ldots, \evali_{n} &\in S_2.
\end{align*}
For example for the spectral clustering algorithm the second to $\dimevec^{\mathrm{th}}$ eigenvectors of the graph shift operators are used to detect communities in networks \cite{Riolo2014} and therefore, $j=1$ is the appropriate offset.
\hfill $\lhd$
\end{requirement}
The DK theorem of Theorem \ref{thm_davis_kahan_Bhatia} is slightly more powerful than the one stated in \cite{vonLuxburg2007}. \cite{vonLuxburg2007} compares the eigenvectors corresponding to eigenvalues of the two matrices which fall within an interval $S_1$. The formulation in Theorem \ref{thm_davis_kahan_Bhatia} allows eigenvalues corresponding to the eigenvectors under comparison to extend beyond $S_1$ as long as they do not enter $S_2;$ the two intervals are separated by an interval of length $\delta.$

\subsection{Bounding the Spaces Spanned by the First $\dimevec$ Eigenvectors} \label{sec_bounding_first_r_evecs}
Using the  DK theorem stated in Theorem \ref{thm_davis_kahan_Bhatia}, we can obtain a bound on the distance of the spaces spanned by the first $\dimevec$ eigenvectors of two symmetric matrices under the following very mild assumption. 

\begin{assumption}\label{ass_nonzero_eigengap}
Assume matrices $\mati, \matii \in \mathbb{R}^{n \times n}$ to both have a non-zero $\dimevec^{\mathrm{th}}$ eigengap, i.e., $\evali_1 \leq \ldots \leq \evali_n$ and $\evalii_1 \leq \ldots \leq \evalii_n$ are such that $\evali_\dimevec \neq \evali_{\dimevec+1}$ and $\evalii_\dimevec \neq \evalii_{\dimevec+1}$.
\hfill $\lhd$
\end{assumption}

In the comparison of the spaces spanned by the first $\dimevec$ eigenvectors it is natural to assume a non-zero $\dimevec^{\mathrm{th}}$ eigengap. If the $\dimevec^{\mathrm{th}}$ and $(\dimevec+1)^{\mathrm{th}}$ eigenvalue are equal then their corresponding eigenvectors are shared. Therefore, an arbitrary choice would have to be made which of the basis elements of the (at least two dimensional) eigenspace corresponding to the $\dimevec^{\mathrm{th}}$ eigenvalue should be considered in the eigenvector comparison and no meaningful comparison could be made.

\begin{remark}
Large eigengaps are commonly used to inform the number of eigenvectors which should be used in graphical analysis \cite{vonLuxburg2007}. Therefore, it is usual to have a large $\dimevec^{\rm th}$ eigengap in the spectra under comparison.
\end{remark}

In Theorem \ref{thm_stucture} we demonstrate that, given Assumption \ref{ass_nonzero_eigengap}, we are always able to compare the spaces spanned by the first $\dimevec$ eigenvectors of two symmetric matrices using the DK theorem.

\begin{theorem}\label{thm_stucture}
Consider the matrices holding the eigenvectors corresponding to the $\dimevec$ smallest eigenvalues of each matrix, namely $\evecsi_0 = [w_1, \ldots, w_\dimevec]\in \mathbb{V}_{n, \dimevec}$ and $\evecsii_0 = [v_1, \ldots, v_\dimevec]\in \mathbb{V}_{n,\dimevec}$. Suppose Assumption \ref{ass_nonzero_eigengap} holds for the matrix spectra under comparison.  
Then, there exists a $Q \in O(\dimevec)$ such that,  
\begin{equation} \label{eqn_cost_structureb}
\left\|\evecsi_0 - \evecsii_0 Q\right\|_F 
\leq c_{n,r}\frac{\left\|\mati - \matii\right\|_2}{\max(\evali_{\dimevec+1} - \evalii_\dimevec, \evalii_{\dimevec+1} - \evali_\dimevec)}.
\end{equation}
\end{theorem}
\begin{IEEEproof}
See Appendix~\ref{app:specialdenom}.
\end{IEEEproof}

\begin{remark}
If $\dimevec=n$, it follows from 
Lemma \ref{lemma_existence_of_Q}, that $\left\|\evecsi_0 - \evecsii_0 Q\right\|_F=0$.
For this choice of $\dimevec$, the bound  in Theorem \ref{thm_stucture} also equals zero since $\min(n, n-\dimevec) = 0.$ When $\dimevec=n$, we have the issue of $\evali_{\dimevec+1}$ and $\evalii_{\dimevec+1}$ not being defined; in \cite[p.~317]{Yu2015} this is addressed by defining $\evali_0 = \evalii_0 = - \infty$ and $\evali_{n+1} = \evalii_{n+1} =  \infty$. In our case, any choice leading to a non-zero denominator in the bound is sufficient to obtain the desired result. 
We hence find that for the degenerate parameter choice, $\dimevec=n$, the bound in Theorem \ref{thm_stucture} is tight.  
\hfill $\lhd$
\end{remark}

In Theorem \ref{thm_stucture} we bound the difference of the spaces spanned by the first $\dimevec$ eigenvectors of any two symmetric matrices, (where the two matrices are required to have a nonzero $\dimevec^{\mathrm{th}}$ eigengap), {\it  using only their spectra}. The denominator of the bound in \eqnref{eqn_cost_structureb} contains a very interesting quantity, namely $\max(\evali_{\dimevec+1} - \evalii_\dimevec, \evalii_{\dimevec+1} - \evali_\dimevec),$ which can be interpreted as a {\it between--matrix eigengap}.

\section{Extension to polynomial mappings and non-zero offsets} \label{sec_eval_map}

\subsection{Polynomial Mappings}
Consider polynomial matrix transformations and their action  on matrix spectra and eigenvectors.

\begin{definition}
\cite[p.~36]{Horn1985} define the evaluation of a polynomial $p(t) = c_l t^l + c_{l-1} t^{l-1} +\ldots + c_1 t + c_0$ at a matrix $\mati$ as
$$
p(\mati) =c_l \mati^l + c_{l-1} \mati^{l-1} +\ldots + c_1 \mati + c_0 I.
$$
\hfill $\lhd$
\end{definition}

\begin{theorem} \label{thm_poly_transf}
\cite[p.~36]{Horn1985} Let $p(\cdot)$ be a given polynomial. If $\evali$ is an eigenvalue of $\mati \in \mathbb{R}^{n \times n}$, while $w$ is an associated eigenvector, then $p(\evali)$ is an eigenvalue of the matrix $p(\mati)$ and $w$ is an eigenvector of $p(\mati)$ associated with $p(\evali)$.
\hfill $\lhd$
\end{theorem}
We are therefore able to transform, for example, the \textit{largest} eigenvalues of any symmetric matrix to be comparable to the \textit{smallest} eigenvalues of another symmetric matrix without altering their corresponding eigenvectors and hence keeping our object of inference, i.e., $\rho_1(\espacei, \espaceii),$ unchanged. 

\subsection{Assumptions on the Spectra}\label{sec_poly_assumptions}

In the proof of Theorem \ref{thm_stucture}, we showed that Assumption \ref{ass_nonzero_eigengap} is sufficient to guarantee the presence of valid DK intervals in the case of the comparison of the spaces spanned by the first $\dimevec$ eigenvectors of two spectra which are ordered and indexed in the same way; this corresponds to a zero offet $(j=0).$ We shall now consider the offset parameter $j$
to be potentially greater than zero, and require non-zero $j^{\mathrm{th}}$ and $(j+\dimevec)^{\mathrm{th}}$ eigengaps of the spectra under comparison.

\begin{assumption} \label{ass_non_zero_jth_eigengaps}
Assume, for given $j\geq 0, \dimevec\geq1$, the eigenvalues $\evali_1 \leq \ldots \leq \evali_n$ of $\mati \in \mathbb{R}^{n\times n}$ and the eigenvalues $\evalii_1 \leq \ldots \leq \evalii_n$ of $\matii \in \mathbb{R}^{n\times n}$ to have a nonzero $j^{\mathrm{th}}$ and $(j+\dimevec)^{\mathrm{th}}$ eigengap, i.e., 
$\evali_{j+1}- \evali_j >0,$ $\evali_{j+\dimevec+1}- \evali_{j+\dimevec} >0,$
$\evalii_{j+1}- \evalii_j >0$ and $\evalii_{j+\dimevec+1}- \evalii_{j+\dimevec} >0.$
\hfill $\lhd$
\end{assumption}

\subsection{Constraints on the Polynomial Mappings}
Since the polynomial transformation can change the ordering of the eigenvalues, i.e., $\evali_i<\evali_j$ does not imply that $p(\evali_i) < p(\evali_j),$ we need to place further assumptions on the transformation parameters to ensure the presence of valid DK intervals. 
Consider a given choice of DK intervals, $S_1 = [a,b], S_2 = \mathbb{R}\backslash (a - \delta, b+\delta)$ for some interval separation $\delta$. We will make use of the following two sets to refer to parts of the transformed spectrum, for given $j\geq 0, \dimevec\geq1$,
\begin{align*}
\mathcal{A}_1 = \Big\{ i \in \{1, \ldots, n\} &\backslash \{j+1, \ldots, j+\dimevec\} \\
&: p(\evali_i)> b \Big\};\\
\mathcal{A}_2 = \Big\{ i \in \{1, \ldots, n\} &\backslash \{j+1, \ldots, j+\dimevec\}\\
&: p(\evali_i) < a \Big\}.
\end{align*}

\begin{constraints} \label{ass_poly_ordering}
For given $j\geq 0, \dimevec\geq1$, let the transformation parameters of $p(\cdot)$ be chosen such that 
\begin{equation}
\mathcal{A}_1 \cup \mathcal{A}_2 = \{1, \ldots, n\} \backslash \{j+1, \ldots, j+\dimevec\}.
\end{equation} 
\hfill $\lhd$
\end{constraints}

There are two possible DK interval choices. $S_1$ can either be defined based on eigenvalues in the spectrum of $p(\mati)$ or $\matii.$ 
The two interval choices result in the following two interval parameter triplets $(a, b, \delta)$, for given $j\geq 0, \dimevec\geq1$,
\begin{align}
  a_1 \!&= \!\!\underset{i \in \{j+1, \ldots, j+\dimevec \}}{\min} \!p(\evali_i), b_1 \!=\!\! \underset{i \in \{j+1, \ldots, j+\dimevec \}}{\max} \!p(\evali_i), \nonumber\\
   \delta_1 &= \min\left(  \evalii_{j+\dimevec+1} - b_1, a_1 - \evalii_j\right) ;\label{eqn_interval_choice1} \\
 a_2 &= \evalii_{j+1},\, b_2 = \evalii_{j+\dimevec}\nonumber\\
 \delta_2 &\!=\! \min\!\!\left[  \underset{i \in \mathcal{A}_1}{\min}~ p(\evali_i)\! -\! b_2, a_2 \!-\!  \underset{i \in \mathcal{A}_2 }{\max}~ p(\evali_i) \right] \label{eqn_interval_choice2} \end{align}

\begin{figure*}[t]
\begin{center}
\includegraphics[scale=\scalevariableS,clip]{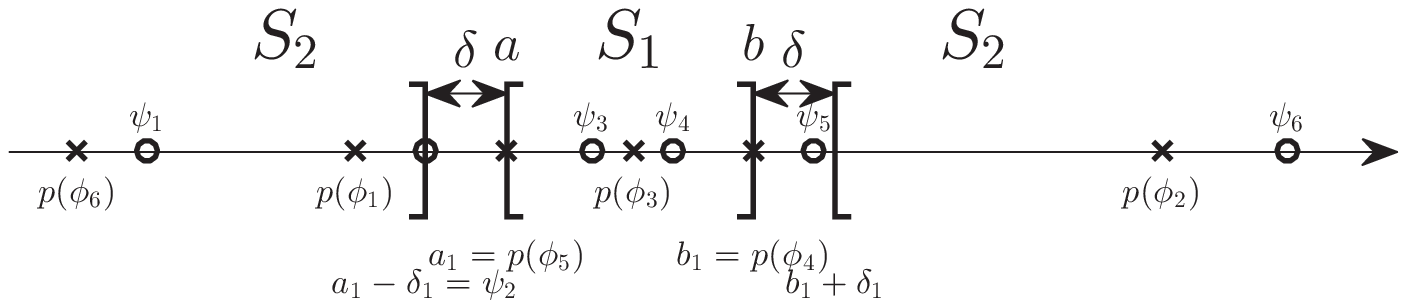}
\caption{A possible eigenvalue configuration used to illustrate Davis--Kahan interval choice \eqnref{eqn_interval_choice1} and the necessary assumptions in order to guarantee its presence.}
\label{fig_choice_of_DK_interval_example_interval_choice_15}
\end{center}
\end{figure*}

\begin{figure*}[t]
\begin{center}
\includegraphics[scale=\scalevariableS,clip]{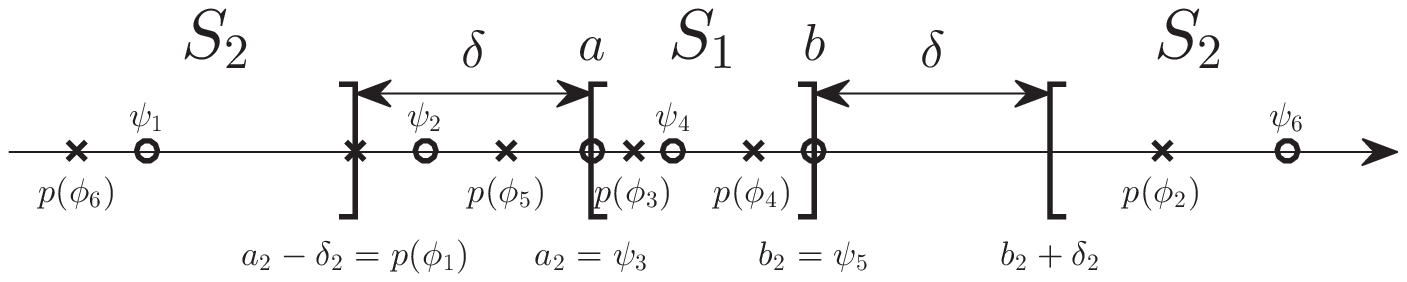}
\caption{A possible eigenvalue configuration used to illustrate the Davis--Kahan interval choice \eqnref{eqn_interval_choice2} and the necessary assumptions in order to guarantee its presence.}
\label{fig_choice_of_DK_interval_example_interval_choice_16}
\end{center}
\end{figure*}

In Fig.~\ref{fig_choice_of_DK_interval_example_interval_choice_15},
we display an example of two spectra composed of 6 eigenvalues each, where the $\evalii_i$'s follow their natural ordering, but the transformed $\evali_i$'s have been given an arbitrary ordering and position for illustration purposes; in practice the choice of $p(\cdot)$ determines the ordering and position.  We choose $\dimevec=3, j=2$, i.e., $\evecsi_2 =[w_3, w_4, w_5]$ and $\evecsii_2 =[v_3, v_4, v_5]$ are the eigenvector matrices to be compared.
For this particular arbitrary set of eigenvalues, interval choice \eqnref{eqn_interval_choice1}, is displayed in Fig.~\ref{fig_choice_of_DK_interval_example_interval_choice_15}, with $a_1 = p(\evali_5), b_1 = p(\evali_4)$ and $\delta_1 = p(\evali_5) -  \evalii_2.$  

In Fig.~\ref{fig_choice_of_DK_interval_example_interval_choice_16} we illustrate interval choice \eqnref{eqn_interval_choice2}, where $a_2 = \evalii_3, b_2 = \evalii_5$ and $\delta_2 = \evalii_3 - p(\evali_1)$. 

Interval choice \eqnref{eqn_interval_choice1} is less favourable than \eqnref{eqn_interval_choice2} since we can clearly observe that $\delta_2>\delta_1$ by comparing 
Figs.~\ref{fig_choice_of_DK_interval_example_interval_choice_15} and \ref{fig_choice_of_DK_interval_example_interval_choice_16}. The purpose of Fig.~\ref{fig_choice_of_DK_interval_example_interval_choice_15} is mainly to demonstrate that both interval choices \eqnref{eqn_interval_choice1} and \eqnref{eqn_interval_choice2} are valid for the displayed set of eigenvalues. For simplicity we will proceed to discuss assumptions only in the context of Fig.~\ref{fig_choice_of_DK_interval_example_interval_choice_16} and interval choice 
\eqnref{eqn_interval_choice2} .

It can be deduced from Fig.~\ref{fig_choice_of_DK_interval_example_interval_choice_16} that the eigenvalues in the illustration satisfy Assumption \ref{ass_non_zero_jth_eigengaps}, i.e., $\evali_{3}- \evali_2 >0,$ $\evali_{6}- \evali_{5} >0,$ $\evalii_{3}- \evalii_2 >0$ and $\evalii_{6}- \evalii_{5} >0.$  (For example, note that $p(\phi_3)\not=p(\phi_2) \implies
\phi_3 > \phi_2,$ because of the eigenvalue ordering, so $\phi_3-\phi_2 >0.$)

Without nonzero $j^{\mathrm{th}}$ and $(j+\dimevec)^{\mathrm{th}}$ eigengaps 
for $\phi_i$ and $\psi_i,$ no valid interval choice could be made, since any violation would immediately lead to a violation of the conditions (Requirements~\ref{rmk_dk_eigenvalue_requirement}) for the presence of valid DK intervals.

For interval  choice \eqnref{eqn_interval_choice2}, Constraints \ref{ass_poly_ordering} says that $\{p(\evali_1), p(\evali_2), p(\evali_6)\}$ must either be larger than $b_2=\evalii_5$ or smaller than $a_2=\evalii_3.$ For the chosen comparison and transformation $p(\cdot),$ we find $\mathcal{A}_1 =\{2\}$ and $\mathcal{A}_2 =\{1, 6\}$ and hence, Constraints \ref{ass_poly_ordering} is satisfied.

As can be seen by observing  $\evalii_5$ in Fig.~\ref{fig_choice_of_DK_interval_example_interval_choice_15} and
$p(\evali_5)$ in Fig.~\ref{fig_choice_of_DK_interval_example_interval_choice_16},
 for interval choices \eqnref{eqn_interval_choice1} and \eqnref{eqn_interval_choice2}, the spectrum not explicitly used in the definition of $S_1$ can extend beyond the boundaries of $S_1$ on both sides. In Section \ref{sec_proof_structure} this issue could be avoided by defining the left boundary of $S_1$ to be $\min(\evali_1, \evalii_1)$; here we  introduce Constraints \ref{ass_poly_non_overlap}. 

\begin{constraints}\label{ass_poly_non_overlap}
$\,$
\begin{enumerate}[label = \Alph*]
\item \label{ass_poly_non_overlap_interval1} For 
interval choice \eqnref{eqn_interval_choice1}, let the transformation parameters of $p(\cdot)$ be chosen such that, for given $j\geq 0, \dimevec\geq1$,
\begin{align*}
\delta_1&>0, \numberthis\label{eqn_delta1_pos}\\
a_1 - \evalii_{j+1} &< \delta_1,\numberthis\label{eqn_use_of_eigengap1} \\
\evalii_{j+\dimevec} - b_1 &< \delta_1.\numberthis\label{eqn_use_of_eigengap2} 
\end{align*}
\item \label{ass_poly_non_overlap_interval2} For interval choice \eqnref{eqn_interval_choice2}, let the transformation parameters of $p(\cdot)$ be chosen such that, for given $j\geq 0, \dimevec\geq1$,
 \begin{align}
\delta_2&>0, \label{eqn_delta2_pos}\\
 a_2 -  \min_{i \in \{j+1, \ldots, j+\dimevec\}} p(\evali_i) &< \delta_2, \label{eqn_overlap_left}\\
 \max_{i \in \{j+1, \ldots, j+\dimevec\}} p(\evali_i) - b_2 &< \delta_2.\label{eqn_overlap_right}
\end{align}
\end{enumerate}
\hfill $\lhd$
\end{constraints}

In essence, in \eqnref{eqn_use_of_eigengap1}, \eqnref{eqn_use_of_eigengap2}, \eqnref{eqn_overlap_left} and \eqnref{eqn_overlap_right} we require that, if any of the eigenvalues with indices $(j+1), \ldots, (j+\dimevec)$ of the spectrum not explicitly used in the definition of $S_1,$ fall outside of $S_1$, then they must not be further than $\delta_i$ away from the boundary of $S_1,$ where $i \in \{1,2\}$ depending on the interval choice we are considering. 
For interval choice \eqnref{eqn_interval_choice2} and the eigenvalues displayed in Fig.~\ref{fig_choice_of_DK_interval_example_interval_choice_16}, Constraints~\ref{ass_poly_non_overlap}\ref{ass_poly_non_overlap_interval2}  requires $\delta_2>0,$ $ \evalii_3 -  p(\evali_5) < \evalii_3 - p(\evali_1)$ and $p(\evali_4) - \evalii_5 < \evalii_3 - p(\evali_1).$

\subsection{The Main Theorem}\label{sec_main_theorem}

Now we are able to extend the DK theorem to be applicable to any two symmetric matrices satisfying Assumption \ref{ass_non_zero_jth_eigengaps}.

\begin{theorem}\label{thm_stucture_transformed}
Let $\mati, \matii \in \mathbb{R}^{n \times n}$ be symmetric matrices  with eigenvalues $\evali_1 \leq \evali_2 \leq \ldots \leq \evali_n$ and $\evalii_1 \leq \evalii_2 \leq \ldots \leq \evalii_n$ and corresponding eigenvectors $\{ w_1, w_2, \ldots, w_n\}$ and $\{ v_1, v_2, \ldots, v_n\}$, respectively. 
Let the matrices holding the eigenvectors corresponding to $\dimevec$ consecutive eigenvalues of each matrix be denoted by $\evecsi_j = [w_{j+1}, \ldots, w_{j+\dimevec}]\in \mathbb{V}_{n,\dimevec}$ and $\evecsii_j = [v_{j+1}, \ldots, v_{j+\dimevec}]\in \mathbb{V}_{n,\dimevec}$. Further, let Assumption \ref{ass_non_zero_jth_eigengaps} hold for the spectra of $\mati$ and $\matii.$ 

Let $p(\cdot)$ be a polynomial transformation satisfying Constraints \ref{ass_poly_ordering} and \ref{ass_poly_non_overlap}\ref{ass_poly_non_overlap_interval1} ($\delta_i =  \delta_1$, 
\eqnref{eqn_interval_choice1})
or Constraints  \ref{ass_poly_ordering} and \ref{ass_poly_non_overlap}\ref{ass_poly_non_overlap_interval2} ($\delta_i = \delta_2$,  
 \eqnref{eqn_interval_choice2})
then, for every unitarily invariant norm, denoted $\|\cdot\|$,
\begin{equation} \label{eqn_cost_structure_transformed_interval1}
\left\|\evecsi_j \evecsi_j^T (I - \evecsii_j \evecsii_j^T)\right\| \leq   \frac{\left\|p(\mati) - \matii\right\|}{\delta_i}. 
\end{equation}

\end{theorem}
\begin{IEEEproof}
See Appendix~\ref{app:mainthm}.
\end{IEEEproof}

We extend our result in Theorem \ref{thm_stucture_transformed} to bound $\left\|\evecsi_j - \evecsii_j Q\right\|_F,$ i.e., $\rho_1,$ in Corollary \ref{cor_main_result}. However, Theorem \ref{thm_stucture_transformed} is equally valid in its own right, if one prefers to work with a bound on the distance of the subspaces in the metric $\left\|\evecsi_j \evecsi_j^T (I - \evecsii_j \evecsii_j^T)\right\|_2,$ i.e., $\rho_2.$

\begin{corollary} \label{cor_main_result}
Let $\mati, \matii, \evecsi_j, \evecsii_j, \delta_i$ and $p(\cdot)$ satisfy the conditions given in Theorem \ref{thm_stucture_transformed}, then, 
there exists a $Q \in O(\dimevec)$ such that
\begin{equation} \label{eqn_cost_structure_transformed}
\left\|\evecsi_j - \evecsii_j Q\right\|_F \leq  c_{n,r}  \frac{\left\|p(\mati) - \matii\right\|_2}{\delta_i}. 
\end{equation}
\end{corollary}
\begin{IEEEproof}
This follows by applying Lemma \ref{lemma_relating_Grassmann_distance_to_DK} to Theorem \ref{thm_stucture_transformed}, noting that the matrix two norm or spectral norm  is unitarily invariant as is required in Theorem \ref{thm_stucture_transformed}.
\end{IEEEproof}

\begin{remark}
For the parameter choice $j =0$ and $p(x) = x$ for $x \in \mathbb{R},$ \eqnref{eqn_cost_structure_transformed_interval1} yields the result in \eqnref{eqn_cost_structureb}, where $\underset{i \in \mathcal{A}_1}{\max}~ p(\evali_i), \evalii_0$ are  defined to equal $- \infty$. Hence, Corollary  \ref{cor_main_result} includes Theorem \ref{thm_stucture} as a special case, i.e., the standard DK theorem is a special case of our extended version, and consequently our bounds are guaranteed to be at least as tight as those given by the standard DK theorem. \hfill $\lhd$
\end{remark}

Analogously to Theorem \ref{thm_stucture}, the bound in \eqnref{eqn_cost_structure_transformed}
{\it  depends only on the spectra} of the two matrices.

Within the restrictions imposed by Constraints \ref{ass_poly_ordering} and  \ref{ass_poly_non_overlap}, we aim to choose the polynomial transformation such that it minimises the bound on the eigenvector difference. We therefore have to choose the polynomial transformation to minimise the numerator of our bound \eqnref{eqn_cost_structure_transformed_interval1}, the spectral norm of the matrix difference $\| p(\mati) - \matii \|_2$, and to maximise the bound's denominator, the maximal eigenvalue interval separation $\delta_1$ or $\delta_2$. 
\smallskip
\noindent\fbox{%
  \parbox{3.38in}{%
    In practice, we recommend choosing  transformations to separately minimise the two subproblems, 
posed by the two possible denominators, and then to work with the smaller of the two bounds in order to achieve an overall minimal bound.
  }%
}

The polynomial transformation of a matrix can also be interpreted as finding the ideal matrix in the sense of producing a minimal bound on the difference of the spaces spanned by the eigenvectors, while preserving the eigenvectors of the untransformed matrix. 

\section{Affine transforms}\label{sec:affinetransforms}
Affine matrix transformations $f(\mati) = c_1 \mati + c_0 I , \,\,(c_1, c_0 \in \mathbb{R})$ are special cases of the polynomial matrix transformations.
We consider the cases $c_1=0, c_1>0$ and $c_1<0$ separately. 
For $c_1=0$ all information in the spectrum is lost, i.e., $ f(\evali_i) = c_1 \evali_i +c_0 =c_0 $, and we cannot find intervals $S_1$ and $S_2$ such that the transformation parameters satisfy Constraints~\ref{ass_poly_ordering}, so there exist no valid DK eigenvalue intervals. 
For $c_1>0$ the ordering of the eigenvalues is preserved in the transformed spectrum, while  for $c_1<0$ the ordering of eigenvalues is reversed in the transformed spectrum. We treat these latter two cases separately throughout the rest of this section. 

For affine transformations, the quantities for the $j^{\mathrm{th}}, (j+1)^{\mathrm{th}}, (j+\dimevec)^{\mathrm{th}}$ and $(j+\dimevec+1)^{\mathrm{th}}$ eigenvalues in the transformed spectrum used in the interval definitions 
 \eqnref{eqn_interval_choice1} and \eqnref{eqn_interval_choice2} 
are displayed in Table~\ref{tab_affine_eigenvalue_simplifications}.
Note that all the quantities in Constraints~\ref{ass_poly_non_overlap} can be found simply by plugging in values from Table~\ref{tab_affine_eigenvalue_simplifications}. 
We see for example, for an affine transformation satisfying Constraints \ref{ass_poly_ordering} and \ref{ass_poly_non_overlap} with  $c_1>0,$ interval choice \eqnref{eqn_interval_choice1} corresponds to the following DK intervals:
 \begin{align*}
 S_1 &= [f(\evali_{j+1}), f(\evali_{j+\dimevec})],\\
 \delta &= \min(\evalii_{j+\dimevec+1} - f(\evali_{j+\dimevec}),  f(\evali_{j+1}) - \evalii_j),\\
 S_2 &= \mathbb{R} \backslash (f(\evali_{j+1}) -\delta_1, f(\evali_{j+\dimevec}) + \delta_1).
 \end{align*}

\begin{table}
\caption{Explicit form of the $j^{\mathrm{th}}, (j+1)^{\mathrm{th}}, (j+\dimevec)^{\mathrm{th}}$ and $(j+\dimevec+1)^{\mathrm{th}}$ eigenvalues in the spectrum of $f(\mati)$, where $f(\evali_i) = c_1 \evali_i +c_0$ is an affine transformation. \label{tab_affine_eigenvalue_simplifications}}
\begin{center}
\[
\begin{array}{|c|| c|c|}
\hline
&c_1>0 & c_1<0\\
\hline
\underset{i \in \mathcal{A}_2 }{\max}~ p(\evali_i) & f(\evali_j) & f(\evali_{j+\dimevec+1})  \\
\underset{i \in \{j+1, \ldots, j+\dimevec\}}{\min} p(\evali_i) & f(\evali_{j+1}) & f(\evali_{j+\dimevec}) \\
\underset{i \in \{j+1, \ldots, j+\dimevec\}}{\max} p(\evali_i) & f(\evali_{j+\dimevec}) &  f(\evali_{j+1}) \\
\underset{i \in \mathcal{A}_1}{\min}~ p(\evali_i) & f(\evali_{j+\dimevec+1})& f(\evali_j) \\
\hline
\end{array}
\]
\end{center}
\end{table}

 Note that the two interval choices \eqnref{eqn_interval_choice1} and \eqnref{eqn_interval_choice2} together with the distinction between affine transformations for $c_1>0$ and $c_1<0$ produce four different possible values for the DK interval separation $\delta,$ which are, from  \eqnref{eqn_interval_choice1}
\begin{align}
\delta_{1,+} &= \min\left(  \evalii_{j+\dimevec+1} - f(\evali_{j+\dimevec}), f(\evali_{j+1}) - \evalii_j\right),\label{eqn_delta_1+}\\
\delta_{1,-} &= \min\left(  \evalii_{j+\dimevec+1} - f(\evali_{j+1}), f(\evali_{j+\dimevec}) - \evalii_j\right),\label{eqn_delta_1-}
\end{align}
and, from  \eqnref{eqn_interval_choice2},
\begin{align}
\delta_{2,+} &= \min\left(  f(\evali_{j+\dimevec+1}) - \evalii_{j+\dimevec}, \evalii_{j+1} - f(\evali_j) \right),\label{eqn_delta_2+}\\
\delta_{2,-} &= \min\left(  f(\evali_j) - \evalii_{j+\dimevec}, \evalii_{j+1} - f(\evali_{j+\dimevec+1}) \right).\label{eqn_delta_2-}
\end{align}

\section{Proof of concept example}\label{sec_d_reg_example}

As a proof of concept we will apply Corollary \ref{cor_main_result} to the comparison of the spaces spanned by eigenvectors corresponding to the 3 smallest eigenvalues of the 
graph shift operators of $d$-regular graphs. Firstly, we briefly define the graph shift operators used in our discussion.

For a general degree matrix $D$, the unnormalised graph Laplacian, $L,$ is defined as 
$
L=D-A,
$
where $A$ is the adjacency matrix. The symmetric 
normalised graph Laplacian is given by 
$L_{sym} = D^{-1/2} L D^{-1/2}.$
When all degrees are equal, i.e., $D =d I$, where $I$ is the identity matrix, we refer to the graph as $d$-regular.

\begin{figure}[t!]
\begin{center}
\includegraphics[scale=0.5,clip]{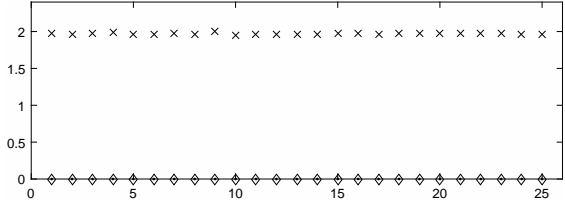}
\caption{ Dots represent attained  distances, $\rho_1(\espacei_0, \espaceii_0),$ diamonds represent bound values calculated from Corollary \ref{cor_main_result}, using an affine matrix transformation, and `x' symbols represent bound values calculated from Theorem \ref{thm_stucture}.}
\label{fig_d_reg}
\end{center}
\end{figure}

In Fig.~\ref{fig_d_reg} we visualise the attained distances, measured in the metric $\rho_1(\espacei_0, \espaceii_0)$ calculated using \eqnref{eqn_stewart_sun_exact_calculation}, of the spaces spanned by the three eigenvectors corresponding to the 3 smallest eigenvalues of 
$L$ and $\lsym$  for 25 30-regular graphs with 300 nodes.
In addition, we show bound values on $\rho_1(\espacei_0, \espaceii_0)$ produced via the standard DK theorem (right hand side of \eqnref{eqn_cost_structureb}) and our theorem with extended applicability (right hand side of \eqnref{eqn_cost_structure_transformed}) using the simplest type of polynomial transformation, namely a linear affine transformation. 

In the case of the graph Laplacian matrices we see that for $d$-regular graphs, the standard DK bound  is non-zero in general. 
Contrariwise, for $d$-regular graphs, $L_{sym}=d^{-1}L,$ i.e., there exists a choice of affine transformation parameters, namely $c_1=d^{-1}, c_0=0,$ 
which maps $L$ exactly to $L_{sym}$, such that 
our extended DK bound is identically zero.

The algorithm discussed in Paper II finds optimal transformation parameters in the case of affine transforms. The result shown in Fig.~\ref{fig_d_reg} was obtained {\it automatically}\/ using this algorithm, without any additional information. The equivalence of the spaces spanned by the first $\dimevec$ eigenvectors was correctly identified by our extended DK bounds methodology.

In the above example both the standard DK theorem and our extended applicability version were able to be applied in comparing $L$ and $\lsym.$ 
In other cases, such as when comparing $A$ and $L,$ the standard DK theorem  cannot even be applied since the required DK intervals of 
Requirements~\ref{rmk_dk_eigenvalue_requirement} do not exist; however, in this case our extended applicability version of the DK theorem can still be used to determine upper bounds. In fact, for $d$-regular graphs, (for which, e.g., $L_{sym} = I- d^{-1} A$), we find that our extended applicability bound identifies the equality of the subspaces spanned by the eigenvectors corresponding to the largest eigenvalues of $A,$ and the spaces spanned by the eigenvectors corresponding to the smallest eigenvalues of the Laplacians $L$ and $\lsym.$

\section{Summary and Conclusions}\label{sec:sum}
In Theorem \ref{thm_stucture} we showed DK bounds always exist for comparison of the spaces spanned by the first $r$ eigenvectors of two symmetric matrices.
Further,  we  proved an extended version of the DK theorem, which applies to the comparison of spaces spanned by any $r$ consecutive corresponding eigenvectors of two symmetric matrices;
Theorem \ref{thm_stucture_transformed} and Corollary \ref{cor_main_result} cover the metrics $\rho_2$ and $\rho_1,$ respectively.

The main tool in the extension of the theorem was the consideration of a polynomial transformation of one of the matrices under comparison. Our extended version includes the original DK theorem as a special case and consequently our bounds are guaranteed to be at least as tight as those given by the original DK theorem. The case of affine transformations as a class of polynomial transformations was discussed. 
A proof of concept example demonstrated that in the case of affine transformations our extended DK theorem, implemented as discussed in detail in Paper II, outperforms the conventional DK theorem by automatically recovering the exact relationship between the adjacency matrix and the unnormalised and normalised graph Laplacians in the case of $d$-regular graphs. 

\section*{Acknowledgment}
The work of Johannes Lutzeyer is supported by the EPSRC (UK).

\appendix{}

\subsection{Proof of Lemma~\ref{lemma_existence_of_Q} } \label{app_existence_of_Q}

Suppose $\evecsi, \evecsii \in \mathbb{V}_{n,\dimevec}$ such that the spaces spanned by their columns, $\espacei, \espaceii \in \mathbb{G}_{n,\dimevec}$, coincide, i.e., $\espacei = \espaceii.$ We want to show that there exists an orthogonal matrix $Q \in O(\dimevec)$ such that $\evecsi = \evecsii Q.$

We begin by extending the bases of $\espacei$ formed by the columns of $\evecsi$ and $\evecsii$ to be bases of $\mathbb{R}^n$ by adding orthogonal columns $\evecsi_\perp, \evecsii_\perp \in \mathbb{V}_{n,n-\dimevec}$ using the Gram-Schmidt algorithm. By joining the columns we obtain the orthogonal matrices $\evecsi_n = [\evecsi, \evecsi_\perp]$ and $\evecsii_n = [\evecsii, \evecsii_\perp]$ .

Consider $Q_n = \evecsii_n^T \evecsi_n$. Note, $\evecsii_n Q_n = \evecsii_n \evecsii_n^T \evecsi_n = \evecsi_n.$ From $\espacei = \espaceii$ it follows that their orthogonal complement spaces coincide, i.e., $\espacei^\perp = \espaceii^\perp$. Since $\evecsi$ and $\evecsii$ span $\espacei$ and $\evecsi_\perp$ and $\evecsii_\perp$ span $\espacei^\perp$, we have by the definition of orthogonal complement spaces that $\evecsi^T \evecsii_\perp  = \evecsii^T \evecsi_\perp  = \mathbf{0}_{\dimevec, n-\dimevec},$ the matrix of all zeros.
Therefore, $Q_n$ is block diagonal,
\begin{align*}
Q_n 
&= \evecsii_n^T \evecsi_n
 = \left[ \begin{smallmatrix} \evecsii^T \\ \evecsii_\perp^T \end{smallmatrix} \right] \left[ \begin{smallmatrix} \evecsi & \evecsi_\perp \end{smallmatrix} \right] 
= \left[ \begin{smallmatrix} \evecsii^T \evecsi & \evecsii^T \evecsi_\perp \\ \evecsii_\perp^T \evecsi & \evecsii_\perp^T \evecsi_\perp \end{smallmatrix} \right]\\
&= \left[ \begin{smallmatrix} \evecsii^T \evecsi & \mathbf{0}_{\dimevec,n-\dimevec} \\ \mathbf{0}_{n-\dimevec,\dimevec} & \evecsii_\perp^T \evecsi_\perp \end{smallmatrix} \right] 
=\left[ \begin{smallmatrix} Q & \mathbf{0}_{\dimevec,n-\dimevec} \\ \mathbf{0}_{n-\dimevec,\dimevec} & Q_\perp\end{smallmatrix} \right] .
\end{align*}

Since $\evecsi_n$ and $\evecsii_n$ are orthogonal, $Q_n$ is also orthogonal. Together the block diagonality and the orthogonality of $Q_n$ imply orthogonality of $Q$. Finally, when only observing the top left $\dimevec \times \dimevec$ block of the equation $\evecsi_n = \evecsii_n Q_n $, we obtain $\evecsii \evecsii^T \evecsi = \evecsi.$ Therefore, we have constructively shown the existence of an orthogonal matrix $Q$ which satisfies $\evecsi = \evecsii Q$.

\subsection{Proof of Lemma~\ref{lemma_relating_Grassmann_distance_to_DK}} \label{app:relatemetrics}
From \eqnref{eqn_defn_subspace_distance} and Proposition 2.2 in \cite[p.~2911]{Vu2013} we have that,  
$$
\frac{1}{2} \left[\inf_{R \in O(\dimevec)} \|\evecsi -\evecsii R \|_F \right]^2 \leq \|\sin\Theta(\espacei, \espaceii)\|^2_F.
$$
Therefore, if we define $Q$ to be the matrix for which the infimum over $O(\dimevec)$ is attained, then, 
\begin{align*}
\!\!\!\!
\left\|\evecsi-\evecsii Q\right\|_F &\leq \sqrt{2} \left\|\sin\Theta(\espacei,\espaceii)\right\|_F \\
&= \sqrt{2} \left\|\evecsi \evecsi^T (I - \evecsii \evecsii^T)\right\|_F \numberthis \label{eqn_vu_projectors}\\
&\leq \!c_{n,r} \left\|\evecsi \evecsi^T (I - \evecsii \evecsii^T) \right\|_2. \numberthis  \label{eqn_Frob_to_spectral_norm}
\end{align*} 
Here \eqnref{eqn_vu_projectors} follows from Equation (2.6) in \cite[p.~2911]{Vu2013} .
\eqnref{eqn_Frob_to_spectral_norm} makes use of the fact that for any matrix $X$ the relation $\left\|X\right\|_F \leq \sqrt{{\rm rank}(X)} \left\|X\right\|_2$ holds \cite[p.~628]{Bernstein2009} and by applying the following simplification: ${\rm rank}\left(\evecsi \evecsi^T (I - \evecsii \evecsii^T)\right) = \min\left({\rm rank}\left(\evecsi \evecsi^T\right),{\rm rank}\left(I - \evecsii \evecsii^T\right)\right) = \min\left(\dimevec, n-\dimevec\right).$

\subsection{Proof of Theorem~\ref{thm_stucture}}\label{app:specialdenom}
By Lemma \ref{lemma_relating_Grassmann_distance_to_DK} there exists a $Q \in O(\dimevec)$ such that,
\begin{equation}\label{eqn_relation_of_norms}
\left\|\evecsi_0 - \evecsii_0 Q\right\|_F \leq  c_{n,r}\left\|\evecsi_0 \evecsi_0^T (I - \evecsii_0 \evecsii_0^T)\right\|_2. 
\end{equation}
Now we bound the term $\left\|\evecsi_0 \evecsi_0^T (I - \evecsii_0 \evecsii_0^T)\right\|_2$ in \eqnref{eqn_relation_of_norms} using the DK Theorem in the form of Theorem \ref{thm_davis_kahan_Bhatia}; in order for the theorem to apply to our comparison of the first $\dimevec$ eigenvectors of two matrices, i.e., for $\evecsi$ to equal $\evecsi_0$ and for $\evecsii$ to equal $\evecsii_0$,  
we need to find intervals such that the conditions identified in 
Requirements~\ref{rmk_dk_eigenvalue_requirement} are satisfied for $j=0$. From the ordering of the eigenvalues it follows that the conditions in Requirements~\ref{rmk_dk_eigenvalue_requirement} simplify to finding intervals $S_1$ and $S_2$ such that the first $\dimevec$ eigenvalues of one of the matrices are contained in $S_1$ -- while the $(\dimevec+1)^{\mathrm{th}}$ eigenvalue is not contained in $S_1$ --  and the last $n-\dimevec$ eigenvalues of the other matrix are contained in $S_2$ -- while the $\dimevec^{\mathrm{th}}$ eigenvalue is not contained in $S_2$.

Either of the two spectra can be used in the definition of $S_1$ and therefore, there exist two valid interval choices with parameters $(a, b, \delta)$ given by,
\begin{align}
a_1 &= \min(\evali_1, \evalii_1), \quad b_1=\evali_\dimevec, \quad
\delta_1 =   \evalii_{\dimevec+1} - \evali_\dimevec; \label{eqn_DK_intervals_1}\\
a_2 &= \min(\evali_1, \evalii_1), \quad b_2=\evalii_\dimevec, \quad
\delta_2 = \evali_{\dimevec+1} - \evalii_\dimevec. 
\label{eqn_DK_intervals_2} 
\end{align}

Parameter choice \eqnref{eqn_DK_intervals_1} leads to interval $S_1 = [\min(\evali_1, \evalii_1), \evali_\dimevec]$ containing the first $\dimevec$ eigenvalues of $\mati$ and not containing $\evali_{\dimevec+1}$ by Assumption \ref{ass_nonzero_eigengap}. Interval $S_2 = \mathbb{R} \backslash (\min(\evali_1, \evalii_1) - \evalii_{\dimevec+1} + \evali_\dimevec, \evalii_{\dimevec+1})$ includes the last $n-\dimevec$ eigenvalues of $\matii$ and excludes $\evalii_\dimevec$ by Assumption \ref{ass_nonzero_eigengap}. 

Similarly, parameter choice \eqnref{eqn_DK_intervals_2} leads to interval $S_1 = [\min(\evali_1, \evalii_1), \evalii_\dimevec]$ containing the first $\dimevec$ eigenvalues of $\matii$ and not containing $\evalii_{\dimevec+1}$ by Assumption \ref{ass_nonzero_eigengap}. Interval $S_2 = \mathbb{R} \backslash (\min(\evali_1, \evalii_1) - \evali_{\dimevec+1} + \evalii_\dimevec, \evali_{\dimevec+1})$ includes the last $n-\dimevec$ eigenvalues of $\mati$ and excludes $\evali_\dimevec$ by Assumption \ref{ass_nonzero_eigengap}.

The parameter choice which results in the greater separation between the two chosen intervals leads to a smaller bound and is hence preferred. By using $\delta = \max(\evali_{\dimevec+1} - \evalii_\dimevec, \evalii_{\dimevec+1} - \evali_\dimevec)$ we ensure that we are working with the interval choice with greater separation.

Assumption \ref{ass_nonzero_eigengap} implies  $\delta = \max(\evali_{\dimevec+1} - \evalii_\dimevec, \evalii_{\dimevec+1} - \evali_\dimevec)>0,$ i.e.,  at least one of the two parameter choices \eqnref{eqn_DK_intervals_1} and \eqnref{eqn_DK_intervals_2} yield valid DK intervals $S_1$ and $S_2$.  To see this, note that Assumption \ref{ass_nonzero_eigengap} guarantees that the following inequalities are strict: $\evali_\dimevec < \evali_{\dimevec+1} $ and $\evalii_\dimevec < \evalii_{\dimevec+1}$. Now we either have $\evali_{\dimevec+1} > \evalii_\dimevec$, which directly implies $\delta>0$ or $\evali_{\dimevec+1} \leq \evalii_\dimevec$, which together with the strict ordering inequalities, implies, $\evali_\dimevec < \evali_{\dimevec+1} \leq \evalii_\dimevec < \evalii_{\dimevec+1}$ and hence, $\evali_\dimevec  < \evalii_{\dimevec+1}$, which means $\delta>0$.

We have demonstrated that we are always able to choose intervals, with separation $\max(\evali_{\dimevec+1} - \evalii_\dimevec, \evalii_{\dimevec+1} - \evali_\dimevec) >0$, satisfying Requirements~\ref{rmk_dk_eigenvalue_requirement}. Hence, continuing from \eqnref{eqn_relation_of_norms}, by the DK theorem and \eqnref{eqn_Bhatia_dk},
\begin{align*}
&\!\!\!\!\!\!c_{n,r}\left\|\evecsi_0 \evecsi_0^T (I - \evecsii_0 \evecsii_0^T)\right\|_2 \\
&\leq c_{n,r}\frac{\left\|\mati - \matii\right\|_2}{\max(\evali_{\dimevec+1} - \evalii_\dimevec, \evalii_{\dimevec+1} - \evali_\dimevec)}.
\end{align*}
so that \eqnref{eqn_cost_structureb} is obtained.

\subsection{Proof of Theorem~\ref{thm_stucture_transformed}}\label{app:mainthm}

We bound the term $\left\|\evecsi_j \evecsi_j^T (I - \evecsii_j \evecsii_j^T)\right\|$  using the DK Theorem  \ref{thm_davis_kahan_Bhatia}. From Theorem \ref{thm_poly_transf} we know that $\mati$ and $p(\mati)$ share eigenvectors and therefore, when bounding$\left\|\evecsi_j \evecsi_j^T (I - \evecsii_j \evecsii_j^T)\right\|$ we are able to apply the DK theorem to the spectrum of $p(\mati)$ instead of the spectrum of $\mati$. 
For the DK theorem to apply we need to show that the interval requirements laid out in Requirements~\ref{rmk_dk_eigenvalue_requirement} are satisfied by interval choices \eqnref{eqn_interval_choice1} and \eqnref{eqn_interval_choice2}. We begin by checking the four conditions laid out in Requirements~\ref{rmk_dk_eigenvalue_requirement} for interval choice \eqnref{eqn_interval_choice1}:

Firstly, the condition $p(\evali_{j+1}), \ldots, p(\evali_{j+\dimevec}) \in S_1$ is always guaranteed by the definition of $S_1= \left[\underset{i \in \{j+1, \ldots, j+\dimevec \}}{\min} p(\evali_i), \underset{i \in \{j+1, \ldots, j+\dimevec \}}{\max} p(\evali_i) \right]$ in interval \eqnref{eqn_interval_choice1}.

Secondly, for interval choice \eqnref{eqn_interval_choice1},  Constraints \ref{ass_poly_ordering} imply that for all $l \in \{1, \ldots, n\} \backslash \{j+1, \ldots, j+\dimevec\}$ we have that either 
\begin{multline} \label{eqn_jth_eigengaps_transformed1}
p(\evali_l)< \underset{i \in \{j+1, \ldots, j+\dimevec \}}{\min} p(\evali_i)\\
 \qquad \text{or} \qquad p(\evali_l)> \underset{i \in \{j+1, \ldots, j+\dimevec \}}{\max} p(\evali_i).
\end{multline}
This implies that 
$$\underset{i \in \mathcal{A}_2 }{\max}~ p(\evali_i) < \min_{i \in \{j+1, \ldots, j+\dimevec\}} p(\evali_i) = a_1$$ and $$b_1 = \max_{i \in \{j+1, \ldots, j+\dimevec\}} p(\evali_i) < \underset{i \in \mathcal{A}_1}{\min}~ p(\evali_i).$$ Hence, Contraints \ref{ass_poly_ordering} guarantee that $p(\evali_{1}), \ldots, p(\evali_{j}), \allowbreak p(\evali_{j+\dimevec+1}), \ldots, p(\evali_{n}) \notin S_1.$

Next, the statements in Constraints \ref{ass_poly_non_overlap}\ref{ass_poly_non_overlap_interval1} can be rearranged to $\delta_1>0,$ $\evalii_{j+\dimevec} < \delta_1 + b_1$ and $ a_1 - \delta_1 < \evalii_{j+1}.$ Since, $S_2 = \mathbb{R}\backslash (a_1- \delta_1, b_1 +\delta_1),$ we find that $\evalii_{j+1}, \ldots, \evalii_{j+\dimevec} \notin S_2$ follows immediately from Constraints \ref{ass_poly_non_overlap}\ref{ass_poly_non_overlap_interval1}.

Finally, from the definition of $\delta_1 = \min\left(  \evalii_{j+\dimevec+1} - b_1, a_1 - \evalii_j\right)$ in interval choice \eqnref{eqn_interval_choice1} it follows that $\delta_1 \leq  a_1 - \evalii_j$,  hence,  $\evalii_j \leq a_1 - \delta_1$ and similarly, $\delta_1\leq \evalii_{j+\dimevec+1} - b_1,$ hence, $b_1 + \delta_1\leq \evalii_{j+\dimevec+1}.$  Therefore, $\evalii_1, \ldots, \evalii_j,\evalii_{j+\dimevec+1}, \ldots,\allowbreak  \evalii_n \in S_2$  follows immediately from the parameter choice \eqnref{eqn_interval_choice1}. 
Therefore, under Constraints  \ref{ass_poly_ordering} and \ref{ass_poly_non_overlap}\ref{ass_poly_non_overlap_interval1} on the transformation parameters, interval choice \eqnref{eqn_interval_choice1} satisfies the conditions laid out in 
Requirements~\ref{rmk_dk_eigenvalue_requirement}. 

In the case of interval choice \eqnref{eqn_interval_choice2}, Assumption \ref{ass_non_zero_jth_eigengaps}, and Constraints \ref{ass_poly_ordering} and \ref{ass_poly_non_overlap}\ref{ass_poly_non_overlap_interval2} suffice to guarantee Requirements~\ref{rmk_dk_eigenvalue_requirement}:

Firstly, the condition $\evalii_{j+1}, \ldots, \evalii_{j+\dimevec} \in S_1$ is always guaranteed by the definition of $S_1= \left[\evalii_{j+1}, \evalii_{j+\dimevec} \right]$ in \eqnref{eqn_interval_choice2}.

Secondly,  the eigengaps in the spectrum of $\matii$ implied by Assumption \ref{ass_non_zero_jth_eigengaps}, 
specifically $\evalii_{j+1}- \evalii_j >0$ and $\evalii_{j+\dimevec+1}- \evalii_{j+\dimevec} >0,$ guarantee that $\evalii_1, \ldots, \evalii_j,\evalii_{j+\dimevec+1}, \ldots,\allowbreak  \evalii_n \notin S_1.$

Next, the statements in Constraints \ref{ass_poly_non_overlap}\ref{ass_poly_non_overlap_interval2} can be arranged to $\delta_2>0,$ $\underset{i \in \{j+1, \ldots, j+\dimevec \}}{\min} p(\evali_i) > a_2 - \delta_2 $ and $ \underset{i \in \{j+1, \ldots, j+\dimevec \}}{\max} p(\evali_i) < b_2 + \delta_2.$ Since, $S_2 = \mathbb{R}\backslash (a_2- \delta_2, b_2 +\delta_2),$ we find that $p(\evali_{j+1}), \ldots, p(\evali_{j+\dimevec}) \notin S_2$ follows immediately from Constraints \ref{ass_poly_non_overlap}\ref{ass_poly_non_overlap_interval1}.

Finally, from the definition of $\delta_2$ in \eqnref{eqn_interval_choice2},
\begin{align*}
\delta_2 &= \min\left(  \underset{i \in \mathcal{A}_1}{\min}~ p(\evali_i) - b_2, a_2 -  \underset{i \in \mathcal{A}_2 }{\max}~ p(\evali_i) \right)  \\
&\leq \underset{i \in \mathcal{A}_1}{\min}~ p(\evali_i) - b_2\\
&\Rightarrow b_2 + \delta_2 \leq \underset{i \in \mathcal{A}_1}{\min}~ p(\evali_i), \numberthis \label{eqn_proof_cond4_1}
\end{align*}
and similarly, 
\begin{align*}
\delta_2 &= \min\left(  \underset{i \in \mathcal{A}_1}{\min}~ p(\evali_i) - b_2, a_2 -  \underset{i \in \mathcal{A}_2 }{\max}~ p(\evali_i) \right)  \\
&\leq  a_2 -  \underset{i \in \mathcal{A}_2 }{\max}~ p(\evali_i)\\
&\Rightarrow a_2 - \delta_2 \geq \underset{i \in \mathcal{A}_2 }{\max}~ p(\evali_i).\numberthis \label{eqn_proof_cond4_2}
\end{align*}
Furthermore, Constraints \ref{ass_poly_ordering} imply that for all $l \in \{1, \ldots, n\} \backslash \{j+1, \ldots, j+\dimevec\}$ either $p(\evali_l) \geq \underset{i \in \mathcal{A}_1}{\min}~ p(\evali_i)$ or $p(\evali_l) \leq\underset{i \in \mathcal{A}_2 }{\max}~ p(\evali_i).$ This together with Equations \eqnref{eqn_proof_cond4_1} and \eqnref{eqn_proof_cond4_2} implies that for all $l \in \{1, \ldots, n\} \backslash \{j+1, \ldots, j+\dimevec\}$ either $p(\evali_l) \geq b_2 + \delta_2$ or $p(\evali_l) \leq a_2 - \delta_2.$ Therefore,  $p(\evali_{1}), \ldots, p(\evali_{j}), \allowbreak p(\evali_{j+\dimevec+1}), \ldots, p(\evali_{n}) \in S_2$  follows from the definition of interval  \eqnref{eqn_interval_choice2} together with Constraints \ref{ass_poly_ordering}.

By Constraints  \ref{ass_poly_non_overlap}, $\delta_1>0$ and $\delta_2>0.$ Therefore, all requirements of the DK theorem are satisfied. Hence,  equation~(\refeq{eqn_cost_structure_transformed_interval1}) holds,
where, for  \eqnref{eqn_interval_choice1}, $\delta_i = \delta_1$, and for \eqnref{eqn_interval_choice2}, $\delta_i = \delta_2.$

\end{document}